\documentclass[11pt]{article}

\makeindex
\usepackage[all]{xy}
\usepackage{graphicx}

\usepackage{amsmath}
\usepackage{amsfonts}
\usepackage{amssymb}
\usepackage{amsthm}
\newcommand{\isom}{\cong}

\newcommand{\Z}{{\bf{Z}}}

\newcommand{\Q}{{\bf{Q}}}
\newcommand{\Qbar}{{\overline{\Q}}}

\newcommand{\C}{{\bf{C}}}
\newcommand{\F}{{\bf{F}}}
\newcommand{\T}{{\bf{T}}}

\newcommand{\m}{{\mathfrak{m}}}

\newcommand{\divs}{\!\mid\!}
\newcommand{\ndiv}{\!\nmid\!}
\newcommand{\tensor}{\otimes}

 1
\DeclareFontEncoding{OT2}{}{} % to enable usage of cyrillic fonts
  
\newcommand{\ord}{{\rm ord}}

\newcommand{\annT}{{{\rm Ann_\T}}}

\usepackage{bm}
\bmdefine\bmu{\mu}
\newcommand{\mur}{{{\bmu}_{r}}}

\newcommand{\comment}[1]{}

\newtheorem{lem}{Lemma}[section]

\newtheorem{prop}[lem]{Proposition}

\newtheorem{thm}[lem]{Theorem}

\theoremstyle{definition}

\newtheorem{rmk}[lem]{Remark}

\newcommand{\thetitle}
{Rational torsion in elliptic curves and the cuspidal subgroup}

\begin{document}
%\ssp
\parindent=2em

\title{\thetitle}
\author{Amod Agashe
\footnote{This material is based upon work supported by the National Science Foundation under Grant No. 0603668.}}
\date{}

\maketitle

%\include{abstract}
%abstract.tex
\begin{abstract}
Let~$A$ be an elliptic
curve over~$\Q$ of square free conductor~$N$. 
We prove that if $A$ has a rational torsion point of prime order~$r$
such that $r$ does not divide~$6N$, then $r$ divides the
order of the cuspidal subgroup of~$J_0(N)$.
\end{abstract}

\section{Introduction}

Let~$A'$ be an elliptic curve over~$\Q$ of square free conductor~$N$
and let $A$ be the optimal curve in the isogeny class (over~$\Q$) of~$A'$.
Let $X_0(N)$ be the modular curve over~$\Q$ associated to~$\Gamma_0(N)$,
and let $J_0(N)$ be its Jacobian. 
By~\cite{breuil-conrad-diamond-taylor}, we may view~$A$ as an abelian
variety quotient over $\Q$ of~$J_0(N)$.
% the modular Jacobian $J_0(N)$, where $N$ is the conductor of~$A$.  
%We assume that the kernel of the map $J_0(N)\to A$ is
%connected, i.e., that~$A$ is an optimal quotient of $J_0(N)$
%(this can always be done by replacing $A$ by an isogenous
%curve if needed). 
By dualizing, 
$A$ can be viewed as an abelian subvariety
of~$J_0(N)$, as we shall do in the rest of this article. 
%Let $X_0(N)$ denote the usual modular curve, whose Jacobian is~$J_0(N)$.
% and $J_0(N)$ be...
%Let $f$ be a newform...and $A$ be the dual of~$A_f$.
The {\em cuspidal subgroup}~$C$ of~$J_0(N)(\C)$
is the group of degree zero divisors on~$X_0(N)(\C)$ that
are supported on the cusps. It is known that $C$ is finite.
%We assume throughout this article that $N$ is a square free integer. 
Since $N$ is square free, the cusps of~$X_0(N)$ are defined over~$\Q$, so
$C \subseteq J_0(N)(\Q)_{\rm tor}$. 

\comment{
Suppose in this paragraph that $N$ is prime.
%When $N$ is prime, 
Then Mazur~\cite{mazur:eisenstein} showed that 
$C = J_0(N)(\Q)_{\rm tor}$, i.e., that the cuspidal
divisors account for all rational torsion in~$J_0(N)$.
He used this to show that if $\beta$ is a prime ideal of~$\T$
lying over an odd rational prime such that $\beta$ is in 
the support of~$J_0(N)(\Q)_{\rm tor}$, then the $\beta$-component
of the Shafarevich-Tate group of~$J_0(N)$ is trivial.
Mazur also related the cuspidal subgroup to the component
group of~$J_0(N)$, and in fact this relation was used
to understand the torsion in~$J_0(N)$. 
%\edit{mention Shah?}
Building on results from~\cite{mazur:eisenstein}, 
Emerton~\cite{emerton:optimal} proved
that 
%if $N$ is prime, then
$A(\Q)_{\rm tor} = A(\Qbar) \cap C$. He also related $A(\Qbar) \cap C$
to the order of the component group of~$A$, and his results
imply a significant cancellation between the orders of the torsion
and component groups of~$A$ in the second part of the Birch
and Swinnerton-Dyer conjecture (see, e.g.~\cite{visfac}).
Thus relations between the cuspidal subgroup and rational torsion
are of great significance. 

We now revert back to our assumption that $N$ is square free
(and not necessarily prime).
}

When $N$
is prime, Mazur~\cite{mazur:eisenstein} showed that 
$C = J_0(N)(\Q)_{\rm tor}$; in particular 
$A(\Q)_{\rm tor} \subseteq C$. 
%i.e., the cuspidal
%divisors account for all rational torsion in~$J_0(N)$
%Using results from~\cite{mazur:eisenstein}, 
%Emerton~\cite{emerton:optimal} proved
%that if $N$ is prime, then
%$A(\Q)_{\rm tor} = A(\Qbar) \cap C$. 
The torsion and cuspidal groups
are of independent interest and importance, and relations between
them are of great significance. For example, using such a relation,
Emerton~\cite{emerton:optimal} showed that when $N$ is prime,
the orders of~$A(\Q)_{\rm tor}$ and the arithmetic component
group of~$A$ are the same, which
implies a significant cancellation 
%between the orders of the torsion and component groups of~$A$ 
in the formula given by the second part of the Birch
and Swinnerton-Dyer conjecture for~$A$ (when $N$ is prime), which is in accord
with the conjecture (see, e.g.~\cite{visfac}).

Based on some numerical
data of Cremona~\cite{cremona:algs}
and Stein~\cite{stein:cuspidal}, we suspect that 
$A(\Q)_{\rm tor} \subseteq C$ more generally when $N$ is
square free, i.e., that again 
the cuspidal
divisors ``explain'' the existence of all the rational torsion points
in~$A$. 
%We have not studied what happens in the non-squarefree~$N$
%case, but one could ask whether in general,
%$A(\Q)_{\rm tor}$ is contained in the cuspidal subgroup (which
%need not be rational). In any case, 
In this paper, we prove 
the following result in this direction: % for squarefree~$N$:

\begin{thm} \label{thm:main}
Recall that~$A'$ is an elliptic curve over~$\Q$ of square free conductor~$N$
and $A$ is the optimal curve in the isogeny class of~$A'$.
%Recall that $N$ is a square free integer.
Suppose $r$ is a prime that does not divide~$6N$. \\
(i) If $r$ divides the order of~$A(\Q)_{\rm tor}$,
then $r$ divides the order of~$A \cap C$ (in particular,
$r$ divides the order of the cuspidal subgroup~$C$). \\
(ii) If $r$ divides the order of~$A'(\Q)_{\rm tor}$, then 
$r$ divides the order of the cuspidal subgroup~$C$.
%Then the prime to~$6N$ part of~$A(\Q)_{\rm tor}$ is contained in~$C$.
\end{thm}

The proof of the theorem is given in Section~\ref{sec:proof}.
The main ingredient in the proof is to show that the hypotheses
imply that the cuspform~$f$ associated to~$A$
is congruent to an Eisenstein series
modulo~$r$. Given such a congruence, and the fact that
$f$ is ordinary at~$r$ (which we show), a result of 
Tang~\cite[Thm~0.4]{tang:congruences}
tells us that $A[r]$ has nontrivial
intersection with a subgroup of the cuspidal group~$C$,
thus giving us statement~(i) of our result. Statement~(ii) follows from
statement~(i) by~\cite[Thm.~1.2]{dummigan:torsion}
which says 
that if $\ell$ is a prime such that
$\ell^2 \ndiv N$ (which holds for $\ell =r$, given our hypothesis), then if $A'$
has a rational torsion point of order~$\ell$, then so does~$A$
(see also Remark~\ref{rmk:dum}).
%While we use
%many techniques from~\cite{mazur:eisenstein}, our methods
%are independent from those of~\cite{emerton:optimal}, even for prime~$N$. 

By\cite[III.5.1]{mazur:eisenstein}, the only primes that can divide 
the order of~$A'(\Q)_{\rm tor}$ are~$2$, $3$, $5$ and~$7$, and moreover
there is a finite list of possibilities for~$A'(\Q)_{\rm tor}$. 
In particular, our theorem gives new information only when
$r$ is~$5$ or~$7$ (and $r \ndiv N$). 
However, we expect
that by doing more work (using ideas from~\cite{mazur:eisenstein}), one
should be able to prove the stronger result that for
every prime~$r \ndiv 6N$, the $r$-primary part of
$A(\Q)_{\rm tor}$ is contained in~$C$; also this result should
hold for higher-dimensional abelian subvarieties~$A$ of~$J_0(N)$ 
associated to newforms. It would also be desirable to see if
the hypothesis that $r \ndiv 6N$ can be removed. 
All this will be the subject
of a future paper. The present article may be viewed
as our first step in relating rational torsion of modular
abelian varieties to the cuspidal subgroup
when $N$ is square free, as well as 
generalizing some of the techniques of Mazur~\cite{mazur:eisenstein}
for prime~$N$ to square free~$N$. 
%Our eventual future goal is to
%generalize appropriately some of the results mentioned
%in the second paragraph of this introduction for prime~$N$
%to square free~$N$.

In any case, the theorem above puts 
restrictions on when $5$ and~$7$ can divide
the order of~$A'(\Q)_{\rm tor}$, and may be useful
in its computations, since
the order of the cuspidal subgroup~$C$ 
can be computed (see, e.g.,~\cite{stein:cuspidal}).
It may also be useful theoretically in certain situations where there
is an explicit formula for the order of~$C$. 
For example, if $N$ is a product
of two primes~$p$ and~$q$, then by~\cite[\S3.4]{chua-ling},
the only odd primes that can divide the order of~$C$ are the 
ones that divide $(p^2-1)(q^2-1)$. As a computational application, 
taking $p = 1013$
and $q = 10007$,
we see that $5$ and~$7$ cannot divide the order of the rational torsion
subgroup of 
any elliptic curve over~$\Q$ of conductor $N = 1013 \cdot 10007$.

The organization of this article is as follows.
In Section~\ref{sec:eis},
we show how to construct certain desirable 
Eisenstein series. In Section~\ref{sec:lemmas}, we
state some 
other preliminary results needed for the proof
of Theorem~\ref{thm:main}. These results are about
certain constraints regarding the Fourier coefficients of~$f$ 
%modulo~$r$
arising out of the existence of rational $r$-torsion, and 
could be of independent interest. Finally,
we give the proof
of Theorem~\ref{thm:main} in Section~\ref{sec:proof}.
Note that in any given section, we continue to 
use the notation introduced in earlier sections. \\
%The proof is given in Section~\ref{sec:proof}. 

%\begin{cor}
%Recall that $N$ is a square free integer. 
%If a prime~$r$ divides the order of~$A(\Q)_{\rm tor}$,
%then $r$ divides $6N \prod_{p \divs N} (p^2 -1)$.
%\end{cor}
%\begin{proof}
%\end{proof}

\noindent {\it Acknowledgement:} We are grateful to Barry Mazur 
for pointing out a construction that we 
used in the proof of Proposition~\ref{prop:eis},
and to Neil Dummigan for conveying the proof of Lemma~\ref{lem:dum},
as well as for some very useful comments on an earlier draft.

\section{Certain Eisenstein series} \label{sec:eis}

%Recall that $N$ is square free. 
If $g = g(z)$ is a modular form,
then we will denote its Fourier expansion $\sum_{n \geq 0} a_n(g) q^n$
at the cusp~$\infty$ (where $q = e^{2 \pi i z}$ as usual)
by $g(q)$. If $n$ is a positive integer, then
$\sigma(n)$ denotes the sum of all the positive divisors of~$n$.

\begin{prop} \label{prop:eis}
Recall that $N$ is square free. 
For every prime~$p$ that divides~$N$,
suppose we are given
an integer $\delta_p \in \{1,p\}$ such that  
$\delta_p = 1$ for at least one~$p$.
Then there is an Eisenstein series~$E$ of weight~$2$ on~$\Gamma_0(N)$ 
which is an eigenfunction for all the Hecke operators such that
for all primes~$\ell \ndiv N$, we have 
$a_\ell(E) = \ell + 1$, and for all primes~$p \divs N$,
we have $a_p(E) = \delta_p$. 
%Moreover, $E$ is an Eisenstein series,
%unless $\delta_p = p$ for all primes~$p$ dividing~$N$.
\end{prop}
\begin{proof}
The normalized Eisenstein series~$e$ of weight~$2$ and level~$1$ has
$q$-expansion $e(q) = 1/24 - \sum_{n \geq 1} \sigma(n) q^n$.
It is not a modular form of level~$1$, but it is an eigenfunction
for all the Hecke operators. We shall construct the desired
Eisenstein series by starting with~$e$ and ``raising the level''.
 
%Let $p$ be a prime dividing~$N$. Then $e(q)$ and~$e(q^p)$
%are both modular forms of level~$p$. They are eigenvectors for
%all Hecke operators, except for~$U_p$. 

Let $g = \sum_{n \geq 0} a_n(g) q^n$ be a normalized eigenfunction 
of some level~$M$
and let $r$ be a prime that does not divide~$M$.
Let $(B_r g)(z) = g(rz)$. Then %at level~$Mr$, on $q$-expansions, 
we have (see, e.g., \cite[p.~141]{atkin-lehner})

\begin{eqnarray} \label{eqn:qexps}
B_r(\sum_{n \geq 0} a_n q^n) & = & \sum_{n \geq 0} a_n q^{nr}, \nonumber\\
U_r(\sum_{n \geq 0} a_n q^n) & = & \sum_{n \geq 0} a_{nr} q^n, \nonumber\\
{\rm and\ } T_\ell(\sum_{n \geq 0} a_n q^n) & = & \sum_{n \geq 0} 
a_{n\ell} q^n +  \sum_{n \geq 0} \ell a_n q^{n \ell},\ \ \ \ \  \forall\ \ell \ndiv Mr,
\end{eqnarray}
where $T_\ell$ and $U_r$ are the usual Hecke operators at level~$Mr$.
For the moment, let $T_r$ denote the $r$-th Hecke operator of level~$M$;
then equation~(\ref{eqn:qexps}) holds for~$T_r$ as well.
Thus from the formulas above, we see that $T_r = U_r + r B_r$.
Since $g$ is an eigenfunction for $T_r$ with eigenvalue $a_r(g)$, 
we deduce that 
$U_r(g) = a_r(g) \cdot  g - r \cdot B_r(g)$ and $U_r(B_r(g))=g$. 
Thus $U_r$ preserves the complex vector
space~$V$ generated by~$g$ and~$B_r(g)$, and
the characteristic
polynomial of $U_r$ on this subspace is 
$U_r^2 - a_r(g) U_r + r$. The elements of~$V$ are
eigenvectors for all the other Hecke operators.
%The other Hecke operators acts
%a multiplication by appropriate scalars on~$V$.
Now suppose $a_r(g) = 1+r$, as will be the case in our application.
Then the 
characteristic polynomial becomes
$U_r^2 - (1+r) U_r + r$, whose 
roots are~$1$ and~$r$.
Thus the action of~$U_r$ is diagonalizable on~$V$. Moreover,
one checks that a basis of normalized eigenvectors
(for all the Hecke operators)
is~$g_r = g - r \cdot B_r(g) = g(q) - r \cdot g(q^r)$ 
and~$\tilde{g}_r  = g - B_r(g) = g(q) - g(q^r)$,
with eigenvalues~$1$ 
and~$r$ respectively for~$U_r$. 
If $g$ is actually a modular form (of level~$M$), then $g_r$ and~$\tilde{g}_r$
are modular forms of level~$Mr$.
Since $g$ is normalized,
$a_r(g_r) = 1$ and $a_r(\tilde{g}_r) = r$. 
Moreover,  
%since raising levels by~$r$ does not disturb the Fourier coefficients at other primes, 
one sees from the construction that for all primes~$\ell \neq r$, we have
$a_\ell(g_r) = a_\ell(\tilde{g}_r) = a_\ell(g)$. 

Now pick a prime~$p$ that divides~$N$ such that $\delta_p = 1$. 
Taking $M = 1$, $r = p$, and $g = e$ in the discussion above, 
and considering that $a_p(e) = \sigma(p) = 1+ p$,
%we see
%that $U_p$ preserves the subspace generated by~$e(q)$ 
%and~$e(q^p) = B_p(e)(q)$,
%and its characteristic polynomial is
%$U_p^2 - a_p(e) U_p + p = U_p^2 - (1+p) U_p + p$.
%The roots of this polynomial are~$1$ and~$p$,
%and thus the action of~$U_p$ is diagonalizable. Moreover,
%one checks that a basis of normalized eigenvectors
%(for all the Hecke operators)
%is~$e_p = e - p \cdot B_p(e) = e(q) - p \cdot e(q^p)$ 
%and~$\tilde{e}_p  = e - B_p(e) = e(q) - e(q^p)$,
%with eigenvalues~$1$ 
%and~$p$ respectively.
we get an Eisenstein series~$e_p$ 
that is an eigenvector for all the Hecke operators
such that 
for all primes~$\ell \neq p$, we have
$a_\ell(e_p) = a_\ell(e) = \sigma(\ell) = \ell + 1$,
and  $a_p(e_p) = 1$.
Moreover, $e_p$ %(q) = e(q) - p e(q^p)$ 
is a modular form of level~$p$
(see, e.g., \cite[p.~47]{diamond-im}).
This proves the desired result if~$N = p$.
Note that unlike~$e_p$, neither~$e$ nor $\tilde{e}_p$ are modular forms,
which is the reason for our hypothesis that 
$\delta_p = 1$ for at least one prime~$p$ dividing~$N$.

If another prime~$s$ divides~$N$, then we apply the procedure
two paragraphs above, taking $M=p$, $r=s$, and 
$g = e_p$. Since $a_s(e_p) = a_s(e) = 1+s$,
we get an eigenform
for all Hecke operators such that 
for all primes~$\ell \ndiv N$, the $\ell$-th Fourier coefficient is
$\ell + 1$, the $p$-th Fourier coefficient is~$1$,
and the $s$-th Fourier coefficient
may be chosen to be~$1$ or~$s$. 
This proves the desired result if~$N = ps$.

If $N$ is a product of more than two distinct primes,
then
%Since $N$ is square free,  and  raising levels by a prime
%that divides~$N$ does not disturb the Fourier coefficients
%at other primes, we see that 
by repeating the procedure
in the previous paragraph for any additional primes that divide~$N$,
%(replacing $q$ by other primes
%dividing~$N$ and using the eigenforms obtained in the
%previous paragraph as starting point) 
we get an eigenform~$E$ with
$a_\ell(E) = \ell +1$ for all primes $\ell \ndiv N$,
$a_p(E)=1$, and for all primes $s \divs N$ with $s \neq p$,
$a_s(E)$ can be chosen to be~$1$ or~$s$.
\end{proof}

The fact that one can construct interesting Eisenstein series 
by raising levels as in the proof above
was pointed out to us by B.~Mazur. In fact, a series as in 
the proposition above was used
for the special case when $N$ is prime in~\cite{mazur:eisenstein} 
(the series~$e'$ 
in~\S~II.5 on p.~78 in loc. cit.). 

\section{Some results on Fourier coefficients} \label{sec:lemmas}
As before, $f$ denotes the cuspform of weight~$2$ on~$\Gamma_0(N)$
associated to~$A$. Then $f$ has integer Fourier coefficients.
Let $w_p$ denote the sign of the Atkin-Lehner involution~$W_p$
acting on~$f$.
In this section, we prove certain results that show how
the existence of rational $r$-torsion
in~$A$ is related to the Fourier coefficients of~$f$.
%modulo~$r$.
\comment{
Some of the discussion
in this section is a bit involved, and the 
reader who is primarily interested in seeing
%only results
%needed from this section to give 
the proof of Theorem~\ref{thm:main} may
jump to Section~\ref{sec:proof} after reading
the statements of Lemma~\ref{lem:mazur0}, Lemma~\ref{lem:dum},
and Proposition~\ref{prop:mazur2}. 
}

The following lemma is perhaps well known.
\begin{lem} \label{lem:mazur0}
Suppose a prime~$r$ divides the order of~$A(\Q)_{\rm tor}$. 
Then for all primes $\ell \ndiv N$, we have
$a_\ell(f) \equiv 1 + \ell \bmod r$ and 
if $p \divs N$, then $a_p(f) = - w_p$.
\end{lem}
\begin{proof}
The proof of the first claim
follows from the discussion in~\cite[p.~112--113]{mazur:eisenstein};
we repeat some of the arguments in loc. cit. for the convenience
of the reader.
Let $P$ be a point of order~$r$ in~$A(\Q)_{\rm tor}$ and 
let $G$ be a finite quotient of~${\rm Gal}(\Qbar/\Q)$
through which the action of~${\rm Gal}(\Qbar/\Q)$ on~$J_0(N)[r]$
factors. Denote by $V$ the $(\T/r\T)[G]$-submodule 
of~$J_0(N)[r]$ generated by~$P$ and by $\m$
the annihilator in~$\T$ of~$V$. Let $S = {\rm Spec\ } \Z$, and let
$J$ denote the N\'eron model of~$J_0(N)$ over~$S$.
Let $V_{/S}$ denote the quasi-finite subgroup scheme 
of~$J[r]$ whose associated 
Galois module is~$V$. Since $N$ is square free, $J_0(N)$ has semi-stable
reduction, and the argument at the bottom of~p.~113 in~\cite{mazur:eisenstein}
shows that $V_{/S}$ is either ${\mathbf \mu}_r  \tensor_{\F_r} \T/\m$
or $\Z_r \tensor_{\F_r} \T/\m$. 
%Considering that $V_{/S}$ has non-trivial
%points over~$\Q$ and $r$ is odd, 
%we see that $V_{/S}$ is $\Z_r \tensor_{\F_r} \T/\m$. 
In either case, 
if $\ell$ is a prime that does not divide~$N$, then 
the Eichler-Shimura
relation 
$T_\ell = {\rm Frob}_\ell + \ell/{\rm Frob}_\ell$ on~$J_{/\F_\ell}$
(where $\ell/{\rm Frob}_\ell$ is the Verschiebung of~$J_{/\F_\ell}$)
%${\rm Frob}_\ell^2 - T_\ell \cdot {\rm Frob}_\ell + \ell = 0$
tells us that 
%\edit{check V irred}
%Let 
%$\ell$ be a prime such that $\ell \ndiv N$.
%Then as shown 
%at the bottom of~p.~113 in~\cite{mazur:eisenstein}
%(where $N$ need not be prime, but must be square free),
%$T_\ell \equiv (1+\ell) \bmod \m$,
%i.e. 
$T_\ell \equiv (1+\ell) \bmod \m$. 
In particular, $T_\ell - (1+\ell)$ annihilates~$P$. 
Since  $T_\ell P = a_\ell(f) P$,
we see that $T_\ell - a_\ell(f)$ annihilates~$P$, and hence
so does
$a_\ell(f) - (1 + \ell)$. But
$P$ has order~$r$, so $r$ divides
$a_\ell(f) - (1 + \ell)$, and hence
$a_\ell(f) \equiv 1 + \ell \bmod r$.

If $p \divs N$, then $a_p(f) = - w_p$ because
$U_p = - W_p$ on the new subspace of~$S_2(\Gamma_0(N),\C)$.
This finishes the proof of the lemma.
%where $W_p$ denotes the Atkin-Lehner involution.
\end{proof}

Keeping in mind the strategy of the proof of 
our main theorem (Theorem~\ref{thm:main})
mentioned 
in the introduction, 
we see from the lemma above and Proposition~\ref{prop:eis}
that coming up with an Eisenstein series~$E$ such that
$a_\ell(f) \equiv a_\ell(E) \bmod r$ for all primes $\ell \ndiv N$
is rather easy. Proving the congruence for all $\ell \divs N$
for a suitable Eisenstein  series
is the tricky part, for which we need the results below.

The following fact is stated without a detailed proof 
in~\cite[\S4]{dummigan:torsion};
the ingredients of the proof were communicated to us by N.~Dummigan.

\begin{lem}[Dummigan] \label{lem:dum}
Let $r$ be an odd prime that divides the order of~$A(\Q)_{\rm tor}$.
If $p$ is a prime that divides~$N$ such that $w_p = 1$, then 
$r \divs (p+1)$.
\end{lem}

\begin{proof}
By the hypothesis, there is a nontrivial  point~$P$ in~$A(\Q)[r]$.
Then $P \in A(\Q_p)[r]$.
Since $p^2 \ndiv N$ (as $N$ is square free) 
and $w_p =1$, the elliptic curve~$A$ has non-split
multiplicative reduction at~$p$. Thus there is a 
$q \in \Q_p^*$ and a Tate
curve~$E_q$ over~$\Q_p$, such that $A$ is isomorphic
to~$E_q$ over an unramified quadratic extension~$K$ of~$\Q_p$.
Now $E_q(\Qbar_p) \isom \Qbar_p/q^\Z$ over~$\Q_p$; let
$x \in \Qbar_p$ be such that its image is in~$E_q(\Qbar_p)[r]$
corresponds to~$P$.
Since $r P =0$, we have $x^r \in q^\Z$, i.e., $x^r =q^n$ for some $n\in \Z$. 
Let $\zeta_r$ be a primitive root of unity in~$\Qbar_p$,
and let $q^{1/r}$ denote a choice of a root of~$X^r = q$ in~$\Qbar_p$.
Then $x = \zeta_r^a q^{b/r}$, for some $a, b \in \{0, \ldots, r-1\}$.

Since $K$ is unramified over~$\Q_p$, its Galois group is generated
by the Frobenius endomorphism, which we will denote by~$\sigma$.
Now $A(\Qbar_p)[r]$ is the same as~$E_q(\Qbar_p)[r]$, 
except that the Galois action on~$A(\Qbar_p)[r]$
is twisted by a nontrivial unramified quadratic character.
%Thus if a non-trivial point~$P \in A(\Qbar_p)[r]$ corresponds to~$x$, then
Thus since $P \in A(\Q_p)$, we have $\sigma(x) = 1/x$ modulo~$q^\Z$.
So the valuation of~$\sigma(x) x$ 
is an integer multiple of that of~$q$,
and since $\sigma$ preserves valuations, 
we have $2 b/r \in \Z$. If $b \neq 0$, then this is possible
only if $r = 2$. 

Now consider the case where $b = 0$. Then $a \neq 0$, and $x = \zeta_r^a$. 
%Then $a \neq 0$, and so $\sigma(\zeta_r $
%then $\zeta_r^r =1 \bmod p$, hence
%$r \divs (p-1)$, and hence $\zeta_r^p = \zeta_r$; in this case,
%since $\sigma$ fixes~$\zeta_p$, we still have 
%$\sigma(\zeta_r) = \zeta_r^p$. 
If $\zeta_r \not\in \Q_p$, then $\sigma(\zeta_r^a) = \zeta_r^{ap}$
and so $\zeta_r^{ap} = 1/\zeta_r^a$.
Since $\zeta_r^a$ is also a primitive~$r$-th root of unity,
we have $r \divs (p+1)$. 
If $\zeta_r \in \Q_p$, then since $\sigma$ fixes~$\zeta_r$, we 
have $\zeta_r = 1/\zeta_r$, i.e., $r=2$.
This proves the lemma.
\end{proof}

\begin{rmk}
In the lemma above, the hypothesis that $r$ is odd is necessary. For example, 
the elliptic curve~$14A1$ has rational $2$-torsion and
$w_2 = 1$ (taking $r = p = 2$, we do not have $r \divs (p+1)$).
\end{rmk}

\comment{
Case II: $b \neq 0$.

Thus $x = $
$\in \mur \cdot \langle q^{1/r} \rangle$ (modulo $q^\Z$),
where $\mur$ is the group of $r$-th roots of unity in~$\Qbar_p$,
and $q^{1/r}$ is a choice of a root of~$X^r = q$ in~$\Qbar_p$.

we
have $E_q[r] \isom (\mur \times \langle q^{1/r} \rangle) q^\Z / q^\Z$.
Since the Frobenius endomorphism~$F$ generates the
Galois group of~$K$ over~$\Q_p$, it acts by inversion
on~$(\mur \times \langle q^{1/r} \rangle)  q^\Z/q^\Z$. Now $F$ acts
by raising to the $p$-th power on~$\mur$, and hence
it fixes an element of $(\mur q^\Z) /q^\Z$ only if $r \divs (p+1)$.
Also, for any $a = 0, \ldots, r - 1$, the valuation of~$F(q^{a/r}) q^{a/r}$ 
is an integer multiple of that of~$q$.
Since $F$ does not change valuations, it fixes some
$q^{a/r}$ only if $2a/r \in \Z$, i.e., $r = 2$.
Thus we see that there is a nontrivial point in 
$A(\Q_p)[r]$ only if $r \divs (p+1)$ or~$r =2$.
Since $A(\Q)[r] \subseteq A(\Q_p)[r]$, we get the statement
in the lemma.
}

Following \cite[p.~77 and p.~70]{mazur:eisenstein},
by a holomorphic modular form in~$\omega^{\tensor 2}$
on~$\Gamma_0(N)$ defined over a ring~$R$,
we mean a modular form in the sense of
\cite[\S1.3]{katz:antwerp350} (see also \cite[\S~VII.3]{deligne-rapoport}).
Thus such an object is a rule which assigns to each pair
$(E_{/T},H)$, where $E$ is an elliptic curve over an $R$-scheme~$T$
and $H$ is a finite flat subgroup scheme of~$E_{/T}$ of order~$N$,
a section of~$\omega_{E/T}^{\tensor 2}$,
where $\omega_{E/T}$ is the sheaf of invariant differentials.
%In our application $R$ will be~$\Z/r\Z$ for some prime $r \ndiv 6N$.
%In this case, 
If $r$ is a prime such that $r \ndiv 6N$ and
$f$ is a modular form of weight~$2$
on~$\Gamma_0(N)$ with coefficients in~$\Z[\frac{1}{6N}]$,
then by~\cite[Lemma~II.4.8]{mazur:eisenstein}, 
there is a holomorphic modular form  in~$\omega^{\tensor 2}$
on~$\Gamma_0(N)$ defined over~$\Z/ r\Z$, which we will
denote $f \bmod r$,  such that the
$q$-expansion of~$f \bmod r$ agrees with the $q$-expansion
of~$f$ modulo~$r$.

\begin{lem}[Mazur] \label{lem:mazur}
Let $R$ be a ring such that~$1/N \in R$. Let $g$
be a holomorphic modular form in~$\omega^{\tensor 2}$
on~$\Gamma_0(N)$ defined over~$R$.
%(notation as in~\cite[p.~77 and p.~70]{mazur:eisenstein}).
Suppose that for some prime~$p$ that divides~$N$,
the $q$-expansion of~$g$ is a power series in~$q^p$,
i.e., there is $h(q) \in R[[q]]$ such that $g(q) = h(q^p)$.
Then $h(q)$ is the $q$-expansion of a 
holomorphic modular form in~$\omega^{\tensor 2}$
on~$\Gamma_0(N/p)$ defined over~$R$.
\end{lem}
\begin{proof}
The lemma is proved in~\cite{mazur:eisenstein} 
under the condition that $N$ is prime, and $p=N$
(Lemma II.5.9 in loc.~cit.).
The same proof works mutatis mutandis to give the lemma above,
with the only change to be made being to replace
certain occurrences of~$N$ by~$p$ (e.g., $q^N$ becomes~$q^p$ everywhere)
and the occurrences of~$N-1$ at the bottom of p.~84 in~\cite{mazur:eisenstein}
by $\phi(N)$, where $\phi$ is the Euler $\phi$-function.
\end{proof}

\begin{prop} \label{prop:mazur2}
Suppose there is a prime~$r$ that does not divide~$6N$
such that $r$ divides the order of~$A(\Q)_{\rm tor}$. 
Then there is a prime~$p$ that divides~$N$ such that $w_p = -1$.
\end{prop}
\begin{proof} 
%Multiply $f$ and~$E$  by $24$ to reduce mod~$r$.
%$In the following, the symbol $\equiv$ will stand for
%``congruent modulo~$r$'', 
%\edit{needed?}
%and the letters $\ell$ and~$p$
%will be reserved for primes. 

Suppose, contrary to the conclusion of the lemma, that for every  prime~$p$ 
that divides~$N$, we have $w_p = 1$. 
If $M$ is a positive integer, then let
us say that a holomorphic modular form~$g$ in~$\omega^{\tensor 2}$
on~$\Gamma_0(M)$ defined over~$\Z/r\Z$ is {\em special at level~$M$}
if $a_n(g) \equiv \sigma(\frac{n}{(n,M)}) \prod_{p \mid M} 
(-1)^{{\rm ord}_p(n)} \bmod r$
for all positive integers~$n$.
Using Lemma~\ref{lem:mazur0} and the fact that $f$ is an eigenvector
for all the Hecke operators,
we see that $f \bmod r$ is special at level~$N$. \\

\noindent{\em Claim:} 
If $M$ is a square free integer and 
$g$ is a holomorphic modular form in~$\omega^{\tensor 2}$
on~$\Gamma_0(M)$ 
defined over~$\Z/r\Z$ that is special at level~$M$ and 
$s$ is a prime that divides~$M$, then there exists
a holomorphic modular form in~$\omega^{\tensor 2}$
on~$\Gamma_0(M/s)$ 
defined over~$\Z/r\Z$ that is special at level~$M/s$
(which is also square free).
\begin{proof}
%Pick a prime~$s$ that divides~$N$.
%We will now construct a 
%a holomorphic modular form~$g$ in~$\omega^{\tensor 2}$
%on~$\Gamma_0(N/s)$ defined over~$\Z/r\Z$ that is special at level~$N/s$.
%Pick one such prime, and call it~$s$.
By Proposition~\ref{prop:eis},
there is an Eisenstein series~$E$ which is an eigenvector
for all the Hecke operators, with $a_\ell(E) = \ell + 1$
for all primes $\ell \ndiv M$,
$a_p(E) = p$ 
for all primes~$p$ that divide~$M$ except $p = s$, and $a_s(E) = 1$. 
%Let $n$ be a positive integer. We can write
%$n = m p_1^{e_1} \cdots p_t^{e_t}$, where
%$p_1 = s, p_2, \ldots, p_t$ are all the prime divisors of~$N$ (for some~$t$),
%the $e_i$'s are nonnegative integers, and $m$ is a positive integer
%coprime to~$N$.
%and let $m$ be a positive integer coprime to~$N$.
%Any positive integer can be written
%as  $m p_i^{e_1} \cdots p_t^{e_t}$.
Let $p_1, \ldots, p_t$ be the distinct primes that divide~$M/s$.
Then for any positive integer~$n$,
$$a_n(E) \equiv \sigma \bigg(\frac{n}{(n,M)} \bigg) 
%\prod_{p \mid M} (1)^{{\rm ord}_s(n)}
\prod_{i=1}^{t}{{p_i}^{{\rm ord}_{p_i}(n)}}
\bmod r.$$
%The $n$-th coefficient of~$f$
%is $\sigma(m) \prod_{i=1}^{t}{(-1)^{e_i}}
% \sigma(m) (-1)^{e_1} \prod_{i=2}^{t}{(-1)^{e_i}}$, while
%that of~$E$ is $\sigma(m) (1)^{e_1} \prod_{i=2}^{t}{(p_i)^{e_i}}$.
Since
by Lemma~\ref{lem:dum}, $p_i \equiv -1 \bmod r$ for $i = 1, \ldots, t$,
%primes~$p$ dividing~$N$)
we see that $a_n( E) \equiv a_n( g) \bmod r$ if $n$ is coprime to~$s$, and thus 
$(E(q)-g(q)) \bmod r$ 
is a power series in $q^s$, i.e., 
there is an $h(q) \in (\Z/r\Z)[[q]]$
with $ h(q^s) \equiv (E(q)-g(q)) \bmod r$.
By Lemma~\ref{lem:mazur}, $h(q)$ is the $q$-expansion
of a 
holomorphic modular form, which we again denote~$h$, in~$\omega^{\tensor 2}$
on~$\Gamma_0(M/s)$ defined over~$\Z/r \Z$.
%Now $a_n(h) \equiv a_{ns}(E - g) \bmod r$.

%If $M$ is a postive integer, then let
%us say that a holomorphic modular form~$g$ in~$\omega^{\tensor 2}$
%on~$\Gamma_0(M)$ is {\em special of level~$M$}
%if 
%$a_n(g) \equiv \sigma(\frac{n}{(n,M)}) \prod_{p \divs M} p^{{\rm ord}_p(n)}$
%for all positive integers~$n$.
%We see that $f$ is special of level~$N$. 
Let $g' = h/2$.
We shall now show that $g'$ is special of level~$M/s$.
Let $n$ be a positive integer, $m' = \frac{n}{(n,s)}$,
and $e = {\rm ord}_s(n)$ (so $n = m' s^e$). Then
% and considering that $g$ is special
%of level~$M$, we see that $a_{m's^{e+1}}(g) = a_{m'}(g) a_{s^{e+1}}(g)$. 
%If ${m'}$ is a positive integer 
%coprime to~$s$, and $e$ is a non-negative integer, then 
\begin{eqnarray}\label{eqn0}
a_n(h)  = 
a_{{m'} s^e}(h)  \equiv  a_{{m'} s^{e+1}}(E - g)  
= a_{{m'} s^{e+1}}(E)  - a_{{m'} s^{e+1}}(g)  \bmod r.
\end{eqnarray}
Now 
$a_n(E) = a_{m'}(E) a_{s^{e+1}}(E)$ since $E$ is an eigenfunction and 
$a_n(g) \equiv a_{m'}(g) a_{s^{e+1}}(g) \bmod r$ since $g$ is special.
Putting this in~(\ref{eqn0}), we get
\begin{eqnarray}\label{eqn1}
a_n(h) & \equiv &
a_{m'}(E) a_{s^{e+1}}(E) - a_{m'}(g) a_{s^{e+1}}(g) \nonumber\\
& \equiv & a_{m'}(g) (a_s(E)^{e+1} - a_s(g)^{e+1})\bmod r,
\end{eqnarray}
where the last congruence follows since $a_{m'}(g) \equiv a_{m'}(E) \bmod r$,
considering that $m'$ is coprime to~$s$. 
Now 
\begin{eqnarray}\label{eqn2}
a_s(E)^{e+1} - a_s(g)^{e+1} = 1-(-1)^{e+1} \equiv 1 -s^{e+1} 
\bmod r, 
\end{eqnarray}
since by Lemma~\ref{lem:dum}, $s \equiv -1 \bmod r$.
Also, 
\begin{eqnarray}\label{eqn3}
1 -s^{e+1} = (1-s) (1+s + \cdots +s^e)
\equiv 2 \sigma(s^e) \bmod r, 
\end{eqnarray}
again considering that
by Lemma~\ref{lem:dum}, $s \equiv -1 \bmod r$.
%$= 2 a_{m'}(g) (1-(-1)^{e+1})/(1 - (-1))
%\equiv 2 a_{m'}(g)  (1 - s^{e+1} )/(1-s) \bmod r$,
%since by Lemma~\ref{lem:dum}, $s \equiv -1 \bmod r$.
%Moreover, $a_{m'}(g)  (1 - s^{e+1} )/(1-s) 
%=a_{m'}(g) (1+ s + s^2 + \cdots + s^e) 
%= a_{m'}(g) \sigma(s^e)$. Thus $a_{{m'} s^e}(h) \equiv 2 a_{m'}(g) \sigma(s^e)
%\equiv 2 \{\sigma(\frac{{m'}}{({m'},M)}) \prod_{p \divs M} p^{{\rm ord}_p(m')} \}
%\sigma(s^e) 
%= 2 \sigma(\frac{{m'}s^e}{({m'}s^e,(M/s))}) \prod_{p \divs (M/s)} 
%p^{{\rm ord}_p({m'}s^e)}$.
Thus putting~(\ref{eqn3}) in~(\ref{eqn2}), and the result in~(\ref{eqn1}),
we get
\begin{eqnarray}\label{eqn4}
a_n(h) \equiv  a_{m'}(g) \cdot 2\sigma(s^e) 
\equiv 2 \sigma\bigg(\frac{m'}{(m',M)}\bigg)  \prod_{p \mid M} (-1)^{{\rm ord}_p(m')}  
\cdot \sigma(s^e) 
\bmod r, 
\end{eqnarray}
where the last congruence follows since $g$ is special at level~$M$.
%Now let $n$ be any Recalling that $g' = h/2$ and that $m'$ is coprime to~$s$, 
%we see that if $n$ is a positive
%integer, then $a_n(g') \equiv 
%\sigma(\frac{n}{(n,M/s)}) \prod_{p \mid (M/s)} (-1)^{{\rm ord}_p(n)}
%\bmod r$,
%i.e., $g'$ is special of level~$M/s$.
%Thus, starting with a form that is special of level~$M$,
%we have constructed a form that is special of level~$N/s$
%for any prime~$s \divs N$.
Now since $n = m' s^e$, with $m'$ coprime to~$s$ and $s \ndiv (M/s)$, we have
\begin{eqnarray}\label{eqn5}
\sigma\bigg(\frac{m'}{(m',M)}\bigg)  \sigma(s^e)  
= \sigma\bigg(\frac{m' s^e}{(m',M)}\bigg) 
= \sigma\bigg(\frac{m' s^e}{(m' s^e,M/s)}\bigg) 
= \sigma\bigg(\frac{n}{(n,M/s)}\bigg) 
\end{eqnarray}
and 
\begin{eqnarray}\label{eqn6}
\prod_{p \mid M} (-1)^{{\rm ord}_p(m')} 
=  \prod_{p \mid M, \ p \neq s} (-1)^{{\rm ord}_p(m's^e)} 
= \prod_{p \mid (M/s)} (-1)^{{\rm ord}_p(n)}.
\end{eqnarray}
Using~(\ref{eqn5}) and~(\ref{eqn6}) in~(\ref{eqn4}), and
recalling that $g' = h/2$, we see that
\begin{eqnarray*}
a_n(g') \equiv 
\sigma\bigg(\frac{n}{(n,M/s)}\bigg) \prod_{p \mid (M/s)} (-1)^{{\rm ord}_p(n)}
\bmod r,
\end{eqnarray*}
i.e., $g'$ is special of level~$M/s$.
\end{proof}

Starting with $f \bmod r$, and repeatedly using the claim,
we see that there is 
a holomorphic modular form that is special of level~$1$, which is nontrivial
since the coefficient of~$q$ is~$1 \bmod r$ for a special form
(of any level). 
%At each stage $f_i$ has non-trivial coefficient of $q$ ($=1$?). 
But by~\cite[Lemma II.5.6(a)]{mazur:eisenstein},
there are no nontrivial holomorphic modular forms
of level~$1$ in $\omega^{\tensor 2}$ defined over 
a field of characteristic other than~$2$ and~$3$.
This contradiction proves the lemma.
\end{proof}

In the proof above, the idea of ``lowering levels'' and getting
a contradiction is taken from  an observation
in~\cite{mazur:eisenstein}, where $N$ is prime and the level
is ``lowered'' only once
(see the proof of Prop.~II.14.1 on p.~114 of loc. cit.).
We noticed that the Fourier coefficients work out so nicely 
(in view of Lemma~\ref{lem:dum}) that
the ``level lowering'' process can be repeated (when $N$ is not necessarily prime),
giving the proof above.

\section{Proof of Theorem~\ref{thm:main}} \label{sec:proof}
%Let $r$ be a prime, and let 
%Suppose an integer~$m$ coprime to~$6N$ 
%divides the order of~$A(\Q)_{\rm tor}$ 
%and 
Recall that the hypotheses are that 
$N$ is a square free integer and $r$ is a prime
such that $r \ndiv 6N$ and $r$ divides
the order of~$A(\Q)_{\rm tor}$. We have
to show that $r$ divides the 
order of the cuspidal subgroup~$C$.

If $p$ is a prime that divides~$N$, then 
let $\delta_p = -w_p$ if $w_p = -1$ and
$\delta_p = p$ if $w_p = 1$.
By  Proposition~\ref{prop:mazur2}, for at least one~$p$,
we have $w_p = -1$, i.e., $\delta_p=1$. 
Hence by Proposition~\ref{prop:eis},
there is an Eisenstein series~$E$ such that
for all primes~$\ell \ndiv N$, we have 
$a_\ell(E) = \ell + 1$, and for all primes~$p \divs N$,
$a_p(E) = 1 = -w_p$ if $w_p = -1$ and 
$a_p(E) = p$ if $w_p = 1$. 
In view of Lemma~\ref{lem:dum}, if $p \mid N$ and $w_p = 1$, 
we have $a_p(E) = p \equiv -1 = -w_p \bmod r$. 

\comment{
$a_\ell(E) = \ell + 1 \equiv a_\ell(f) \bmod r$, and for all primes~$p \divs N$,
$a_p(E) = 1 = -w_p \equiv a_p(f) \bmod r$ if $w_p = -1$ and 
$a_p(E) = p$ if $w_p = 1$.
In view of Lemma~\ref{lem:dum}, if $p \mid N$ and $w_p = 1$, 
we have $a_p(E) = p \equiv -1 = a_p(f) \bmod r$. 
$a_n(f) \equiv a_n(E) \bmod r$ for all $n \geq 1$.
}

Considering that $f$ and~$E$ are eigenfunctions for all
the Hecke operators, we see from the paragraph above 
and by Lemma~\ref{lem:mazur0}
that $a_n(f) \equiv a_n(E) \bmod r$ for all $n \geq 1$.
Hence $(f(q) - E(q))\bmod r$ is a constant; call this constant~$c$.
Since $r \ndiv 6N$, %by~\cite[Lemma~II.4.8]{mazur:eisenstein}, 
we may consider 
the holomorphic modular form $(f-E) \bmod r$ in~$\omega^{\tensor 2}$
on~$\Gamma_0(N)$ defined over~$\Z/r \Z$.
Using 
Lemma~\ref{lem:mazur}, for any prime~$p$ dividing~$N$
%we can ``lower the level'' of~$(f-E) \bmod r$ by~$p$, to 
we get a holomorphic modular form in~$\omega^{\tensor 2}$
on~$\Gamma_0(N/p)$ defined over~$\Z/r \Z$, whose $q$-expansion
is the same constant~$c$.
By repeating this process (which we can do since
at each stage we have a $q$-expansion
that is constant -- in fact, the same constant~$c$), we get
%see that $(f(q) -E(q))\bmod r$ is the $q$-expansion of a 
a holomorphic modular form in~$\omega^{\tensor 2}$
on~$\Gamma_0(1)$ defined over~$\Z/r\Z$, whose $q$-expansion
is~$c$. 
By~\cite[Lemma II.5.6(a)]{mazur:eisenstein},
there are no nontrivial holomorphic modular forms
of level~$1$ in $\omega^{\tensor 2}$ defined over 
a field of characteristic other than~$2$ and~$3$.
%But this is not possible
%by~\cite[Lemma II.5.6(a)]{mazur:eisenstein}, considering that $r \neq 2, 3$.
Thus $c \equiv 0 \bmod r$, and so
$a_n(f) \equiv a_n(E) \bmod r$ for $n=0$ as well.
Hence $f \equiv E \bmod r$.

To~$E$ is associated 
a subgroup~$C_E$ of~$C$ by Stevens (see~\cite[Def.~1.8.5]{stevens:thesis}
and~\cite[Def.~4.1]{stevens:cuspidal}).
Since $r \ndiv N$, by Lemma~\ref{lem:mazur0},
$a_r \equiv (1+r) \equiv 1 \bmod r$; in particular,
$f$ is ordinary at~$r$.
By~\cite[Thm~0.4]{tang:congruences}, $A[r] \cap C_E \neq 0$,
and thus $r$ divides the order of~$A \cap C$.
The fact that  $r$ divides the order of~$C_E$ 
follows from the intermediate result
Prop.~1.9 of~\cite{tang:congruences} as well. 
This proves part~(i) of Theorem~\ref{thm:main}.
As mentioned in the introduction, part~(ii) follows from part~(i) 
by taking $\ell =r$ in~\cite[Thm.~1.2]{dummigan:torsion} 
(Dummigan's theorem in turn follows from
the proof of Prop.~5.3 in~\cite{vatsal:mult}).

\begin{rmk} \label{rmk:dum}
Neil Dummigan remarked to us that one need not 
use~\cite[Thm.~1.2]{dummigan:torsion} to deduce part~(ii)
of Theorem~\ref{thm:main} from part~(i) since our methods
prove the following special case of~\cite[Thm.~1.2]{dummigan:torsion}:
if $A'$ is an elliptic curve of square free conductor~$N$ having
a rational point of order~$r$ for a prime~$r$ 
such that $r \ndiv 6N$, then the optimal curve~$A$ in the
isogeny class of~$A'$ also has a rational point of order~$r$. 
Clearly, this shows that part~(i) implies part~(ii).
The details of what we claimed two sentences above are as follows: 
Lemma~\ref{lem:mazur0} holds with $A$ replaced by $A'$
under the additional hypothesis that $r$ is odd
(by considering reduction modulo~$\ell$, we see that 
if $\ell$ is a prime such that $\ell \ndiv N$, then $r$ divides
$|A'(\F_\ell)| = |A(\F_\ell)| = a_\ell(f) - (1 + \ell)$) and
Lemma~\ref{lem:dum} also holds with $A$ replaced by $A'$
(the hypothesis $w_p=1$ implies that $A$ has non-split multiplicative
reduction; hence so does~$A'$ and the proof goes through
with $A$ replaced by $A'$). In the proofs of
Proposition~\ref{prop:mazur2} and Theorem~\ref{thm:main},
the only place where the hypothesis that~$A$ has a rational
point of order~$r$ is used is in quoting Lemmas~\ref{lem:mazur0}
and~\ref{lem:dum}. Since the conclusions of these Lemmas hold
under the hypothesis that $A'$ (instead of~$A$)
has a rational point of order~$r$
(and the hypothesis that $r$ is odd, which is already assumed
in Theorem~\ref{thm:main}), the proof of Theorem~\ref{thm:main}
goes through to prove that $A$ has a rational point of order~$r$,
as claimed.
\end{rmk}

\comment{
Suppose an integer~$m$ coprime to~$6N$ 
divides the order of~$A(\Q)_{\rm tor}$ 
and let $P$ be a point of order~$m$. 
It suffice to show that $P \in C$.

Let $G$ be some finite quotient of~${\rm Gal}(\Qbar/\Q)$
through which the action of~${\rm Gal}(\Qbar/\Q)$ on~$J_0(N)[m]$
factors. Let $V$ be the $\T/m\T[G]$-submodule 
of~$J_0(N)[m]$ generated by~$P$. 
Then as shown in Mazur,
for all $\ell \ndiv N$, $T_\ell - (1+\ell)$ annihilates~$V$,
and in particular~$P$; Since $T_\ell P = a_\ell P$, we 
deduce that $a_\ell \equiv 1 + \ell \bmod m$.
Moreover, if $p \divs N$, then $a_p = - w_p$.
By Proposition~\ref{prop:eis}, taking $\delta_p = -w_p$,
there is an Eisenstein series~$E$ such that
$f \equiv E \bmod m$. 

To~$E$ is associated 
a subgroup~$C_E$ of~$C$ by Stevens.
By Tang, for any prime~$r$, $T_r$ acts on~$C_E$ by
multiplication by~$a_r(E)$. 
Let $\m$ be a max id...need to define Eis ideal, 
and check its kernel, in order to take care of
squares of max ids, etc.
}

%\end{proof}

\bibliographystyle{amsalpha}         

\providecommand{\bysame}{\leavevmode\hbox to3em{\hrulefill}\thinspace}
\providecommand{\MR}{\relax\ifhmode\unskip\space\fi MR }
\providecommand{\MRhref}[2]{%
  \href{http://www.ams.org/mathscinet-getitem?mr=#1}{#2}
}
\providecommand{\href}[2]{#2}

\end{document}